
  \newcount\fontset
  \fontset=1
  \def\dualfont#1#2#3{\font#1=\ifnum\fontset=1 #2\else#3\fi}

  \dualfont\bbfive{bbm5}{cmbx5}
  \dualfont\bbseven{bbm7}{cmbx7}
  \dualfont\bbten{bbm10}{cmbx10}
  \font \eightbf = cmbx8
  \font \eighti = cmmi8 \skewchar \eighti = '177
  \font \eightit = cmti8
  \font \eightrm = cmr8
  \font \eightsl = cmsl8
  \font \eightsy = cmsy8 \skewchar \eightsy = '60
  \font \eighttt = cmtt8 \hyphenchar\eighttt = -1
  
  \font \sixbf = cmbx6
  \font \sixi = cmmi6 \skewchar \sixi = '177
  \font \sixrm = cmr6
  \font \sixsy = cmsy6 \skewchar \sixsy = '60
  \font \tensc = cmcsc10
  
  \font \titlefont = cmr7 scaled \magstep4
  \scriptfont \bffam = \bbseven
  \scriptscriptfont \bffam = \bbfive
  \textfont \bffam = \bbten

  \font\rs=rsfs10 

  \newskip \ttglue

  \def \eightpoint {\def \rm {\fam0 \eightrm }%
  \textfont0 = \eightrm
  \scriptfont0 = \sixrm \scriptscriptfont0 = \fiverm
  \textfont1 = \eighti
  \scriptfont1 = \sixi \scriptscriptfont1 = \fivei
  \textfont2 = \eightsy
  \scriptfont2 = \sixsy \scriptscriptfont2 = \fivesy
  \textfont3 = \tenex
  \scriptfont3 = \tenex \scriptscriptfont3 = \tenex
  \def \it {\fam \itfam \eightit }%
  \textfont \itfam = \eightit
  \def \sl {\fam \slfam \eightsl }%
  \textfont \slfam = \eightsl
  \def \bf {\fam \bffam \eightbf }%
  \textfont \bffam = \eightbf
  \scriptfont \bffam = \sixbf
  \scriptscriptfont \bffam = \fivebf
  \def \tt {\fam \ttfam \eighttt }%
  \textfont \ttfam = \eighttt
  \tt \ttglue = .5em plus.25em minus.15em
  \normalbaselineskip = 9pt
  \def \MF {{\manual opqr}\-{\manual stuq}}%
  \let \sc = \sixrm
  \let \big = \eightbig
  \setbox \strutbox = \hbox {\vrule height7pt depth2pt width0pt}%
  \normalbaselines \rm }


  \def \Headlines #1#2{\nopagenumbers
    \advance \voffset by 2\baselineskip
    \advance \vsize by -\voffset
    \headline {\ifnum \pageno = 1 \hfil
    \else \ifodd \pageno \tensc \hfil \lcase {#1} \hfil \folio
    \else \tensc \folio \hfil \lcase {#2} \hfil
    \fi \fi }}

  \def \Title #1{\vbox{\baselineskip 20pt \titlefont \noindent #1}}

  \def \Date #1 {\footnote {}{\eightit Date: #1.}}

  \def \Authors #1{\bigskip \bigskip \noindent #1}

  \long \def \Addresses #1{\begingroup \eightpoint \parindent0pt
\medskip #1\par \par \endgroup }

  \long \def \Abstract #1{\begingroup \eightpoint
  \bigskip \bigskip \noindent
  {\sc ABSTRACT.} #1\par \par \endgroup }


  \def \lcase #1{\edef \auxvar {\lowercase {#1}}\auxvar }

  \def \goodbreak {\vskip0pt plus.1\vsize \penalty -250 \vskip0pt
plus-.1\vsize }

  \newcount \secno \secno = 0
  \newcount \stno

  \def \seqnumbering {\global \advance \stno by 1
    \number \secno .\number \stno }

  \def\section #1{\global\def\SectionName{#1}\stno = 0 \global
\advance \secno by 1 \bigskip \bigskip \goodbreak \noindent {\bf
\number \secno .\enspace #1.}\medskip \noindent \ignorespaces}

  \long \def \sysstate #1#2#3{\medbreak \noindent {\bf \seqnumbering
.\enspace #1.\enspace }{#2#3\vskip 0pt}\medbreak }
  \def \state #1 #2\par {\sysstate {#1}{\sl }{#2}}
  \def \definition #1\par {\sysstate {Definition}{\rm }{#1}}
  \def \remark #1\par {\sysstate {Remark}{\rm }{#1}}


  \def \proof {\medbreak \noindent {\it Proof.\enspace }}
  \def \proofend {\ifmmode \eqno \square \else \hfill \square
\looseness = -1 \medbreak \fi }

  \def \$#1{#1 $$$$ #1}
  \def\=#1{\buildrel #1 \over =}

  \def\iItem {\smallskip}
  \def\Item #1{\smallskip \item {#1}}
  \newcount \zitemno \zitemno = 0
  \def\izitem {\zitemno = 0}
  \def\zitem {\global \advance \zitemno by 1 \Item {{\rm(\romannumeral
\zitemno)}}}

  \newcount \footno \footno = 1
  \newcount \halffootno \footno = 1
  \def\footcntr {\global \advance \footno by 1
  \halffootno =\footno
  \divide \halffootno by 2
  $^{\number\halffootno}$}
  \def\fn#1{\footnote{\footcntr}{\eightpoint#1}}


  \def \N {{\bf N}}
  \def \C {{\bf C}}
  \def \R {{\bf R}}
  
  \def \({\left (\vrule height 9pt width 0pt}
  \def \){\right )}
  \def \[{\left \Vert }
  \def \]{\right \Vert }
  \def \*{\otimes }
  \def \+{\oplus }
  \def \:{\colon }
  \def \<{\left \langle }   \def \<{\langle }
  \def \>{\right \rangle }  \def \>{\rangle }
  \def \text #1{{\rm #1}}
  
  \def \curly#1{\hbox{\rs #1\/}}
  \def \ds{\displaystyle}
  \def \and {\hbox {,\quad and \quad }}
  
  \def \calcat #1{\,{\vrule height8pt depth4pt}_{\,#1}}

  \def \for #1{,\quad \forall\,#1}
  \def \inv {^{-1}}
  
  \def \square {\hbox {$\sqcap \!\!\!\!\sqcup $}}
  \def \stress #1{{\it #1}\/}

  \def \|{\Vert }
  \def \inv {^{-1}}



  \def \ifn #1{\expandafter \ifx \csname #1\endcsname \relax }
  \def \cite #1{\ifn{showref}{\rm [\bf #1\rm ]}\else {\tt [\#\string #1]}\fi }
  \def \label #1{\global \edef #1{\number \secno \ifnum \number \stno
= 0\else .\number \stno \fi }\ifn{showlabel}\else {\tt [*\string #1] }\fi }

  \def \scite #1#2{\cite {#1{\rm \hskip 0.7pt:\hskip 2pt #2}}}
  \def \lcite #1{\ifn{showcit}(#1)\else (#1 {\tt\string #1})\fi }

  \newcount \bibno \bibno =0
  \def \newbib #1{\global \advance \bibno by 1 \edef #1{\number \bibno
}}
  \def \bibitem #1#2#3#4{\smallskip \item {[#1]} #2, ``#3'', #4.}
  \def \references {
    \begingroup 
    \bigskip \bigskip \goodbreak 
    \eightpoint 
    \centerline {\tensc References}
    \nobreak \medskip \frenchspacing }

  \def\titletext{AF-algebras and the tail-equivalence relation on
Bratteli diagrams}\footnote{\null}
  {\eightrm 2000 \eightsl Mathematics Subject Classification:
  \eightrm 
  46L05, 
  46L85. 
  }

  \def\Bratteli{B}
  
  \def\Renault{R}
  \def\EL{EL}
  \def\Watatani{W}


  \def\C{C(\Omega)}
  \def\CG{C^*(\G)}
  \def\D{{\cal D}}
  \def\Vertices{{\cal V}}
  \def\Edges{{\cal E}}
  \def\nl{\hfill\cr\cr}

  \def\Ex{\curly {E}}

  \def\G{{\cal G}}
  \def\Real{{\bf R}}
  \def\Complex{{\bf C}}
  
  \def\R{\curly{R}}
  \def\Subalg#1{C(\Omega;\R_{#1})}
  \def\a{\alpha}
  \def\b{\beta}
  \def\d{\delta}
  \def\g{\gamma}
  \def\preE{E^0}
  \def\n#1{\raise 1pt \hbox{$\scriptscriptstyle \#$}#1}
  \def\et{{\hat e}}
  \def\ecan{e}
  \def\eb{{\check e}}
  \def\Kt{{\hat {\cal K}}}
  \def\Kb{{\check {\cal K}}}
  \def\Toep{\curly{T}(\R,\Ex)}
  \def\CRE{C^*(\R,\Ex)}
  \def\bool#1{[{\scriptstyle #1}]\ }

  \Headlines
  {\titletext}
  {R.~Exel and J.~Renault}

  \Title{\titletext}

  \Date{04 May 2003}

  \Authors
  {R.~Exel\footnote{*}{\eightrm Partially supported by CNPq.}
  and
  J.~Renault}

  \Addresses
  {Departamento de Matem\'atica,
  Universidade Federal de Santa Catarina,
  Florian\'opolis,
  Brasil
  (exel@mtm.ufsc.br),
  \par
  D\'epartement de Math\'ematiques,
  Universit\'e d'Orl\'eans,
  France
  (renault@labomath.univ-orleans.fr).}

  \Abstract {Given a Bratteli diagram $\D$ we consider the compact
topological space $\Omega$ formed by all infinite paths on $\D$.   Two such
path are said to be tail-equivalent when they ``have the same tail'',
i.e. when they eventually coincide.   This equivalence relation is
approximately proper and hence one may consider the C*-algebra
associated to it according to a procedure recently introduced by the
first named author and A.~Lopes.  The main result of this work is the
proof that this algebra is isomorphic to the AF-algebra associated to
the given Bratteli diagram.}

  \section{Introduction and statement of the main result}
  Let $\D = (\Vertices,\Edges)$ be a directed graph, where $\Vertices$ is the set of
vertices and $\Edges$ the set of edges.
  Recall (see \cite{\Bratteli}) that $\D$ is a Bratteli diagram if:
  \iItem
  \Item {(a)} one is given a decomposition
  $$
  \Vertices=\bigcup_{n\in\N}\Vertices_n
  $$
  (we adopt the convention according to which the set $\N$ of natural
numbers starts with zero) of $\Vertices$ as the union of pairwise disjoint,
finite, nonempty sets $\Vertices_n$,
  \Item {(b)} for every edge $\varepsilon\in\Edges$, if the source $s(\varepsilon)$ of $\varepsilon$
lies in $\Vertices_n$ then its range $r(\varepsilon)$ lies in $\Vertices_{n+1}$,
  \Item {(c)} for every $n$, the set of edges from $\Vertices_n$ to
$\Vertices_{n+1}$ is finite.
  \medskip We will moreover assume, for simplicity, that:

  \iItem
  \Item {(d)} $\Vertices_0$ is a singleton,
  \Item {(e)} every vertex is the source of some edge,
  \Item {(f)} every vertex is the range of some edge, except for the
single vertex in $\Vertices_0$.

\medskip Assumptions (e) and (f) above correspond to the fact that 
the embeddings between the finite dimensional sub-algebras of the
associated AF-algebra \cite{\Bratteli} are unital and injective.

For each $n\in\N$, let us denote by $\Edges_n$ the subset of $\Edges$ given by
  $$
  \Edges_n = \big\{\varepsilon\in\Edges : s(\varepsilon) \in \Vertices_n\big\}.
  $$
  By (b) above one obviously has that $r(\varepsilon)\in\Vertices_{n+1}$ for every
$\varepsilon\in\Edges_n$.  Therefore $\Edges_n$ is precisely the finite set referred to
in (c).


By a \stress{path} in $\D$ we will mean, as usual, a finite or
infinite sequence $\a=(\a_0,\a_1,\a_2,\ldots)$ of edges such that the
range of each edge $\a_n$ coincides with the source of the following
edge in the sequence.  For a finite path 
  $\a = (\a_0,\a_1,\a_2,\ldots,\a_n)$
  we will say that the \stress{length} of $\a$ is the number of edges involved
(as opposed to the number of vertices). 

If a path $\a$ is such that $s(\a_0)\in\Vertices_0$ then it clearly follows
that $\a_n\in\Edges_n$ for every $n$.  If moreover $\a$ is infinite
then it may be considered as an element of the Cartesian
product
  $$
  \prod_{n\in\N} \Edges_n.
  $$
  Equipping each $\Edges_n$ with the discrete topology it becomes a
compact space by (c) and hence the Cartesian product above is compact
for the product topology.  It is not hard to see that the subset
  $$
  \Omega =
  \big\{\a \in \prod_{n\in\N} \Edges_n: \a \hbox{ is a path}\big\}
  $$
  is a closed subset and hence also compact.

The main object of interest in this work is the equivalence relation
on $\Omega$, sometimes called the \stress{tail-equivalence relation},
defined by
  $$
  \a \sim \b
  $$
  if and only if there exists $n\in\N$ such that $\a_k=\b_k$, for all
$k\geq n$.

Quite commonly one finds that the quotient topological space
$\Omega/{\sim}$ is very badly behaved, often being chaotic (the only
open sets being the empty set and the whole space).  This is the case
e.g.~for the Bratteli diagram of the CAR algebra, namely the diagram
in which each $\Vertices_n$ consists of a single vertex and each
$\Edges_n$ consists of exactly two edges (necessarily both joining the
vertex in $\Vertices_n$ to the vertex in $\Vertices_{n+1}$).

Returning to the general case consider, for each $n\in\N$, the
equivalence relation on $\Omega$ defined by
  $$
  \a \sim_n \b
  $$
  if and only if $\a_k=\b_k$, for all $k\geq n$.

It is apparent that the equivalence class of each $\a\in\Omega$ under
``$\sim_n$'' is determined by the infinite sub-path
$(\a_n,\a_{n+1},\ldots)$ and that the quotient space $\Omega/\sim_n$
is homeomorphic to the space of all infinite paths starting at some
vertex in $\Vertices_n$.
  So this quotient is a well behaved Hausdorff space.  In other words
``$\sim_n$'' is a proper equivalence relation.

Denote by
  $$
  \R = \big\{(\a,\b)\in \Omega\times\Omega: \a\sim\b\big\}
  \and
  \R_n = \big\{(\a,\b)\in \Omega\times\Omega: \a\sim_n\b\big\}.
  $$
  Observe that, according to the strictly technical definition of equivalence
relations, $\R$ \underbar{is} the equivalence relation ``$\sim$'' and
$\R_n$ \underbar{is} ``$\sim_n$''.  From now on we will therefore
refer to ``$\sim$'' and ``$\sim_n$'' as $\R$ and $\R_n$, respectively.

  Notice that $\R$ is the increasing union of the $\R_n$ and hence
$\R$ is an approximately proper equivalence relation.

Since the $\R_n$ are proper it is easy to see that the sub-C*-algebra
  $$
  \Subalg{n}
  $$
  of $\C$ formed by the functions which are constant on
each $\R_n$-equivalence class is *-isomorphic to $C(\Omega/\R_n)$.
  Our next immediate goal will be to describe a certain conditional
expectation from $\C$ onto $\Subalg{n}$.

  For each $n\in\N$ and each $\a\in\Omega$ we will let
  $$
  \R_n(\a) = \big\{ \b\in\Omega: \b\sim_n\a\big\},
  $$
  that is, $\R_n(\a)$ is the equivalence class of $\a$ relative to $\R_n$.
  Observe that a path $\b$ lies in $\R_n(\a)$ if and only if it is of
the form
  $$
  \b = (\g_0,\g_1,\ldots,\g_{n-1},\a_n,\a_{n+1},\ldots),
  $$
  where $\g = (\g_0,\g_1,\ldots,\g_{n-1})$ is any path joining (the single
vertex in) $\Vertices_0$ to the source of $\a_n$.

  For each vertex $v\in\Vertices$\/ let 
  $
  \n{v} 
  $
  denote the number of paths joining $\Vertices_0$ to $v$.  Clearly
$\n{v}$ is finite and nonzero for every $v$. It is also clear from our
discussion above that the number of elements in $\R_n(\a)$ coincides
with $\n{s(\a_n)}$.

\state \label \EContinuous Proposition Given $n\in\N$ and $f\in
\C$, the complex valued function $\preE_n(f)$ defined on
$\Omega$ by
  $$
  \preE_n(f)\calcat \a =
  \sum_{\b\in\R_n(\a)} f(\b)
  \for \a\in\Omega,
  \eqno {(\dagger)}
  $$
  is continuous.

\proof Suppose first that $f(\a)$ depends only on the first $m$
coordinates of $\a$, that is
  $$
  f(\a) = g(\a_0,\a_1,\a_2,\ldots,\a_m)
  \for \a\in\Omega,
  $$
  where $g$ is some complex function defined on $\prod_{k=0}^m \Edges_k$.
Supposing without loss of generality that $m>n$, one easily sees that
$\preE_n(f)\calcat \a$ likewise depends only on the first $m$
coordinates of $\a$.  So $\preE_n(f)$ is continuous.

Returning to the general case observe that the sup-norm
  $
  \|\preE_n(f)\|_\infty
  $
  is bounded by $K\|f\|_\infty$, where $K$ is the maximum number of summands in
\lcite{$\dagger$} (that is, 
  $K= \max\{\n{v}: v\in \Vertices_n\}$).  Now apply the Stone--Weierstrass
Theorem to write $f$ as the uniform limit of a sequence $\{f_k\}_k$
formed by functions $f_k$ each of which depends only on finitely many
coordinates of its argument.  We conclude that $\preE_n(f)$ is the
uniform limit of the sequence $\{\preE_n(f_k)\}_k$ and hence that
$\preE_n(f)$ is continuous.
  \proofend

\state Proposition Given $n\in\N$ and $f\in \C$ consider the
complex valued function $E_n(f)$ defined on $\Omega$ by
  $$
  E_n(f)\calcat \a =
  {1 \over \n{s(\a_n)}} \sum_{\b\in\R_n(\a)} f(\b).
  \eqno {(\ddagger)}
  $$
  Then $E_n$ is a conditional expectation from $\C$ onto
$\Subalg{n}$.  Moreover for every $m\geq n$ one has that $E_nE_m = E_mE_n =
E_m$.  In particular $E_n$ commutes with $E_m$.

\proof 
Notice that $\n{s(\a_n)}$ is a continuous function of $\a$
since it depends only on one  coordinate of $\a$, namely  $\a_n$.  It
then follows from \lcite{\EContinuous}  that
  $E_n(f) = \n{s(\a_n)}\inv \preE_n(f) \in \C$.

Observe that the number of summands in \lcite{$\ddagger$} is exactly
$\n{s(\a_n)}$.  Therefore $E_n(f)\calcat\a$ is precisely the
arithmetic mean of the values of $f$ on the equivalence class of $\a$
relative to $\R_n$.  So it easily follows that $E_n$ is a conditional
expectation onto $\Subalg{n}$.

Given that $n\leq m$ it is clear that $\Subalg{n} \supseteq \Subalg{m}$ and hence
$E_n$ coincides with the identity on $\Subalg{m}$.  It follows that
$E_nE_m=E_m$.  

It remains to prove that $E_mE_n=E_m$.  In order to do so let
$\Vertices_n=\{v_1,\ldots,v_p\}$ and for every $i=1,\ldots,p$ let
$X_i$ denote the set of all paths from $\Vertices_0$ to $v_i$.
  Given $\a\in\Omega$,
  let $Y_i$ be the set of all paths from $v_i$ to $s(\a_m)$.
Therefore a path $\b$ is in $\R_m(\a)$ if and only if $\b$ is of the
form
  $$
  \b = xy\bar\a
  \quad \hbox{(juxtaposition of paths),}
  $$
  where, for some $i=1,\ldots,p$, one has that $x\in X_i$, $y\in Y_i$,
and $\bar\a=(\a_n,\a_{n+1},\ldots)$.

  Fix $i=1,\ldots,p$, and $y\in Y_i$.  We claim that for
any $f\in \C$, we have that
  $$
  \sum_{x\in X_i} f(xy\bar\a)=
  \sum_{x\in X_i} E_n(f)\calcat{xy\bar\a}.
  \eqno{(\star)}
  $$
  In order to see why this is so choose  $x_0\in X_i$.   The
$\R_n$-equivalence class of $x_0y\bar\a$ is therefore formed by all paths
of the form $xy\bar\a$, for $x\in X_i$.
Therefore we have by definition that
  $$
  E_n(f)\calcat{x_0y\bar\a} =
  {1 \over \n{v_i}} \sum_{x\in X_i} f(xy\bar\a).
  $$
Replacing $f$ by $E_n(f)$ above, and observing that $E_n^2=E_n$, we have
  $$
  E_n(f)\calcat{x_0y\bar\a} =
  {1 \over \n{v_i}} \sum_{x\in X_i} E_n(f)\calcat{xy\bar\a},
  $$
  and $(\star)$ follows by applying the definition of $E_n(f)$ on the
left hand side above.

Returning to the proof that  $E_m=E_mE_n$ observe that 
  $$
  E_m(f)\calcat \a =
  {1\over \n{s(\a_m)}} \sum_{\b\in\R_m(\a)} f(\b) =
  {1\over \n{s(\a_m)}} 
    \sum_{i=1}^p \sum_{y\in Y_i} \sum_{x\in X_i} f(xy\bar \a) =$$$$
  \={(\star)}
  {1\over \n{s(\a_m)}} 
    \sum_{i=1}^p \sum_{y\in Y_i} \sum_{x\in X_i}
E_n(f)\calcat{xy\bar\a} =
  E_m(E_n(f))\calcat \a.
  \proofend
  $$


It is our goal in this work to describe the C*-algebra associated to
$\R$ and the collection of conditional expectations $\Ex = \{E_n\}_n$
under the procedure described in \cite{\EL} and to prove it to be
isomorphic to the AF-algebra arising from the Bratteli diagram $\D$.

Recall from \cite{\EL} that the Toeplitz algebra of the pair
$(\R,\Ex)$,  denoted $\Toep$, is the universal C*-algebra
generated by a copy of $\C$ and a sequence $\{\et_n\}_{n\in\N}$ of
projections  subject to the relations:
  \izitem
  \zitem
  $\et_0=1$,
  \zitem $\et_{n+1}\et_n = \et_{n+1}$,
  \zitem
  $\et_n f \et_n = E_n(f) \et_n,$
  \medskip \noindent
  for all $f\in \C$ and $n\in\N$.  

\smallskip 
As in \cite{\EL}, for each $n\in\N$ we will denote by
$\Kt_n$ the closed linear span of the set
  $$
  \{f\et_n g: f,g\in\C\}
  $$
  within $\Toep$.  By \scite{\EL}{2.4} we have that $\Kt_n$ is a
*-subalgebra of $\Toep$.

Still according to \cite{\EL} a \stress{redundancy} is, by
definition, a finite sequence
  $$
  (k_0,\ldots,k_n)\in \prod_{i=0}^n\Kt_i,
  $$
  such that $\sum_{i=0}^n k_ix=0$, for all $x\in \Kt_n$.

  The ideal of $\Toep$ generated by the sums $\sum_{i=0}^n k_i$, for
all redundancies $(k_0,\ldots,k_n)$, is called the \stress{redundancy
ideal}.  The C*-algebra for the pair $(\R,\Ex)$, denoted $\CRE$, is
then defined \scite{\EL}{2.7} as the quotient of $\Toep$ by the
redundancy ideal.  We will denote by $\ecan_n$ the image of $\et_n$ in
$\CRE$.

The main goal of this work is therefore to prove:

  \state Theorem
  \label\Main
  Let $\D$ be a Bratteli diagram and let $\R$ be the tail-equivalence
relation on the infinite path space of\/ $\D$.  If $\Ex = \{E_n\}_n$ is
the sequence of conditional expectations defined above then $\CRE$ is
isomorphic to the AF-algebra associated to the Bratteli diagram $\D$.

\section{Systems of matrix units in $\CRE$}
  In this section we will introduce systems of  matrix units within
$\CRE$ which will correspond to the matrix units of the AF-algebra
associated to the Bratteli diagram $\D$ and will eventually lead us to
the proof of \lcite{\Main}.

For each finite path
  $\g = (\g_0,\g_1,\g_2,\ldots,\g_{n-1})$ starting at $\Vertices_0$ 
and each $\a\in\Omega$ we will say that
  $$
  \g\leq\a
  $$
  when
  $(\a_0,\a_1,\a_2,\ldots,\a_{n-1}) = \g$, that is,
  whenever $\a$ ``starts with $\g$''.
  Moreover we will
  let $I_\g$ be the characteristic function of the set
  $$
  \Omega_\g = \big\{\a\in\Omega: \g\leq \a \big\}.
  $$
  We may then write
  $$
  I_\g(\a) = \bool{\g\leq\a}\!\!
  \for \a\in\Omega,
  $$
  where the brackets stand for the obvious boolean valued function.

For each vertex $v\in\Vertices_n$ let $I^v$ be the characteristic
function
of the  set
  $$
  \Omega^v =   \big\{\a\in\Omega: s(\a_n) = v\big\}.
  $$
  Both $\Omega_\g$ and $\Omega^v$ are clopen sets and hence the
corresponding characteristic functions $I_\g$ and $I^v$ are
continuous.

Before we proceed it will be convenient to extend the notion of the 
``range of an edge'' in order to apply it to finite paths also.  So, 
given a finite path $\g = (\g_0,\g_1,\g_2,\ldots,\g_{n-1})$ we will
let
  $$
  r(\g) := r(\g_{n-1}),
  $$
  so that  the range of a finite path is understood to be the range of its
last edge, as one would naturally expect.  

If $\g = (\g_0,\g_1,\g_2,\ldots,\g_{n-1})$ is a finite path starting
at $\Vertices_0$ and $v\in\Vertices_n$ it is easy to see that
  $$
  I_\g I^v =
  \left\{\matrix{
  I_\g, & \hbox{ if } r(\g) = v, \cr\cr
  0, & \hbox{ if } r(\g) \neq v.
  }\right.
  $$
  We will write this as 
  $$
  I_\g I^v = 
  \bool{r(\g)=v} I_\g.
  \eqno{(\seqnumbering)}
  \label\IgIv
  $$

  \state Lemma 
  \label \EnIg
  Let $\g = (\g_0,\g_1,\g_2,\ldots,\g_{n-1})$ be a finite path
starting at $\Vertices_0$.  Then
  $$
  E_n(I_\g) = {I^{r(\g)} \over \n{r(\g)}}.
  $$

  \proof
  This is just a computation based on the definition of $E_n$ and is
left for the reader.
  \proofend

We are now ready to introduce the matrix units which constitute the
core of this section.

  \definition
  Given paths 
  $\g = (\g_0,\g_1,\g_2,\ldots,\g_{n-1})$  and
  $\d = (\d_0,\d_1,\d_2,\ldots,\d_{n-1})$  of the same length $n$ let 
  $e^n_{\g,\d}$ be the element of $\CRE$ given by
  $$
  e^n_{\g,\d}=
  \n{r(\g)} I_\g \ecan_n I_\d,
  $$
  where $\ecan_n$ is the image of $\et_n$ in $\CRE$, as already
mentioned.

  \state Lemma 
  \label\ZeroEs
  Let $\g$, $\d$, and $n$ be as above.  Then $e^n_{\g,\d}=0$ whenever
$r(\g) \neq r(\d)$.

  \proof
  It suffices to show that $(I_\g \ecan_n I_\d)(I_\g \ecan_n I_\d)^* =0$.  In
order to prove this notice that
  $$
  (I_\g \ecan_n I_\d)(I_\g \ecan_n I_\d)^* =
  I_\g \ecan_n I_\d \ecan_n I_\g \$=
  I_\g E_n(I_\d) \ecan_n I_\g \={(\EnIg)}
  {1 \over \n{r(\d)}} I_\g I^{r(\d)} \ecan_n I_\g \={(\IgIv)}
  {1 \over \n{r(\d)}} \bool{r(\g) = r(\d)} I_\g \ecan_n I_\g.
  \proofend
  $$

  One should therefore concentrate on the $e^n_{\g,\d}$ for
which $r(\g)=r(\d)$.  For these we have:

  \state Lemma 
  \label\MatrixFormulas
  Let $\g$, $\d$, $\zeta$, $\eta$ be finite paths of length $n$ with
  $r(\g) = r(\d)$ and   $r(\zeta) = r(\eta)$.  Then
  $$
  e^n_{\g,\d} e^n_{\zeta,\eta} = \bool{\d=\zeta}   e^n_{\g,\eta}.
  $$

  \proof
  We have
  $$
  \matrix{e^n_{\g,\d} e^n_{\zeta,\eta} &=& \n{r(\g)} \n{r(\zeta)}I_\g \ecan_n I_\d I_\zeta \ecan_n I_\eta \nl
  & = & \bool{\d=\zeta} \n{r(\g)}^2 I_\g \ecan_n I_\d \ecan_n I_\eta \nl
  & = & \bool{\d=\zeta} \n{r(\g)}^2 I_\g E_n(I_\d) \ecan_n I_\eta \nl
  & = & \bool{\d=\zeta} \n{r(\g)}^2 I_\g \ds {I^{r(\d)}\over \n{r(\d)}} \ecan_n I_\eta \nl
  & = & \bool{\d=\zeta} \n{r(\g)} I_\g I^{r(\d)} \ecan_n I_\eta \nl
  & = & \bool{\d=\zeta} \n{r(\g)} I_\g \ecan_n I_\eta \nl
  & = & \bool{\d=\zeta}   e^n_{\g,\eta},\nl
  }
  $$
  concluding the proof.
  \proofend

  \section{Finite Index}
  There is not much more one can say about the present situation
before proving that the expectations $E_n$ introduced above are of
index-finite type according to \cite{\Watatani}.  This must therefore
be our next goal.

\state Proposition \label\QB
Let $n\in\N$ be fixed and let $\Omega_n$ denote the set of all finite paths
  of length $n$  starting at $\Vertices_0$.  For each $\g\in\Omega_n$ let
  $$
  u_\g = \sqrt{\n{r(\g)}}\ I_\g.
  $$
  Then the set $\{u_\g\}_{\g\in\Omega_n}$ is a quasi-basis for
$E_n$.  In particular $E_n$ is of index-finite type.

  \proof
  Let $\g\in\Omega_n$ and let $f\in\C$.  Then, for every
$\a\in\Omega$, we have
  $$
  \matrix{I_\g E_n(I_\g f)\calcat \a & = &
  \bool{\g\leq\a} \ds{1\over   \n{s(\a_n)}} 
    \sum_{\b\in\R_n(\a)} \bool{\g\leq\b}f(\b) \nl
  &=&  \bool{\g\leq\a} \ds{1\over   \n{r(\g)}}
    \ f(\a).\nl
  }
  $$
  It follows that
  $$
  \sum_{\g\in\Omega_n} u_\g E_n(u_\g f)\calcat \a =
  \sum_{\g\in\Omega_n} \n{r(\g)}I_\g E_n(I_\g f)\calcat \a =
  \sum_{\g\in\Omega_n} \bool{\g\leq\a} f(\a) = f(\a),
  $$
  so that   
  $\ds\sum_{\g\in\Omega_n} u_\g E_n(u_\g f) = f$, as desired.
  \proofend

One of the first consequences of this is given in our next:

\state Proposition \label \IandE For every $\g\in\Omega_n$ one has
that $I_\g = e^n_{\g,\g}$.

\proof Let   $\{u_\d\}_{\d\in\Omega_n}$ be the quasi-basis for $E_n$
given by \lcite{\QB}.  By \scite{\EL}{6.2.i} we have that
  $$
  1 = \sum_{\d\in\Omega_n} u_\d \ecan_n u_\d.
  $$
  It follows that
  $$
  I_\g =
  \sum_{\d\in\Omega_n} I_\g u_\d \ecan_n u_\d =
  u_\g \ecan_n u_\g =
  \n{r(\g)}I_\g \ecan_n I_\g = 
  e^n_{\g,\g}.
  \proofend
  $$

By \scite{\EL}{3.7} we have that the canonical embedding of $A$ into
$\CRE$ is injective.  Therefore, since  $I_\g$ is obviously non-zero,
it follows that 
  $$
  e^n_{\g,\d}\,
  e^n_{\d,\g} = 
  e^n_{\g,\g} =
  I_\g \neq 0,
  $$
  whenever $\g,\d\in\Omega_n$ are such that $r(\g)=r(\d)$.  In
particular   $e^n_{\g,\d}\neq0$.

  Fixing $v\in\Vertices_n$ the sub-C*-algebra of $\CRE$ generated by
the set
  $\{e^n_{\g,\d}: r(\g)=r(\d)=v\}$
  is therefore easily seen to be isomorphic to
  $
  M_{\n{v}}(\Complex)
  $
  by \lcite{\MatrixFormulas}.
  If $\Vertices_n=\{v_1,\ldots,v_p\}$ then the
sub-C*-algebra generated by
  $\{e^n_{\g,\d}: \g,\d\in\Omega_n\}$ is isomorphic to
  $$
  M_{\n{v_1}}(\Complex) \oplus M_{\n{v_2}}(\Complex) \oplus \cdots
\oplus M_{\n{v_p}}(\Complex),
  \eqno{(\seqnumbering)}
  \label\SumMatrix
  $$
  by \lcite{\ZeroEs}.

In order to see how the $e^n_{\g,\d}$ relate to each other for
different $n$'s we need a quasi-basis for the restriction of $E_{n+1}$ to
$\Subalg n$.  By \scite{\EL}{6.1} such a quasi-basis is simply given by
$\{E_n(u_\g)\}_{\g\in\Omega_{n+1}}$.  In order to obtain a
concrete expression for the 
$E_n(u_\g)$ we need to introduce yet another characteristic function
of interest.  Given an edge $\varepsilon\in\Edges_n$ we will denote by 
  \def\vI#1{\null^{#1}\kern-2ptI}%
  $\vI{\varepsilon}$
  the
characteristic function of the clopen set
  $$
  \big\{ \a\in\Omega: \a_n = \varepsilon \big\}.
  $$
  Thus
  $$
  \vI{\varepsilon}(\a) = \bool{\a_n=\varepsilon}
  \for \a\in\Omega.
  $$
  Observing that $\vI{\varepsilon}(\a)$  depends only on $\a_n$, and hence is
constant on the $\R_n$-equivalence class of $\a$, it is clear that $\vI{\varepsilon}\in\Subalg n$.

  Given $\g\in\Omega_{n+1}$, say
  $\g= (\g_0,\g_1,\ldots,\g_{n-1},\g_n)$, observe that
  $$
  I_\g = I_{\g'}\,\vI{\g_n},
  $$
  where $\g'=(\g_0,\g_1,\ldots,\g_{n-1})$.  So
  $$
  E_n(I_\g) = 
  E_n(I_{\g'}\,\vI{\g_n}) = 
  E_n(I_{\g'})\,\vI{\g_n} \={(\EnIg)}
  {I^{r(\g')} \over \n{r(\g')}}\,\vI{\g_n} = 
  {\vI{\g_n} \over \n{r(\g')}}.
  $$
  It follows that
  $$
  E_n(u_\g) =
  \sqrt{\n{r(\g)}}\ E_n(I_\g)= 
  {\sqrt{\n{r(\g)}} \over \n{r(\g')}}\ \vI{\g_n}.
  $$
  By \scite{\EL}{6.2.i} we conclude that
  $$
  \ecan_n =
  \sum_{\g\in\Omega_{n+1}} E_n(u_\g) \ecan_{n+1} E_n(u_\g) =
  \sum_{\g\in\Omega_{n+1}} {\n{r(\g)} \over \n{r(\g')}^2}\, \vI{\g_n}\ \ecan_{n+1}\ \vI{\g_n}.
  $$
  If $\zeta,\eta\in\Omega_n$ are such that $r(\zeta)=r(\eta)$
we then have that
  $$
  e^n_{\zeta,\eta}=
  \n{r(\zeta)} I_\zeta \ecan_n I_\eta= 
  \sum_{\g\in\Omega_{n+1}} 
    {\n{r(\zeta)}\n{r(\g)} \over \n{r(\g')}^2}\,
    I_\zeta\vI{\g_n}\ \ecan_{n+1}\ I_\eta \vI{\g_n} \$=
  \sum_{\g\in\Omega_{n+1}} 
    {\n{r(\zeta)}\n{r(\g)} \over \n{r(\g')}^2}
    \bool{r(\zeta)=s(\g_n)}
    I_{\zeta\g_n}\ecan_{n+1}I_{\eta\g_n} 
= \cdots
  $$
  Denoting by $\Edges_{\zeta}$ the set of all edges whose source is
$r(\zeta)$, so that the juxtaposition $\zeta\varepsilon$ lies in $\Omega_{n+1}$ 
if and only if $\varepsilon\in\Edges_{\zeta}$, the above equals
  $$
  \cdots =
  \sum_{\varepsilon\in \Edges_{\zeta}} \ \
  \sum_{\buildrel {\scriptstyle\g\in\Omega_{n+1}} \over 
        {\vrule height 7pt width 0pt \g_n=\varepsilon}} 
    {\n{r(\zeta)}\n{r(\g)} \over \n{r(\g')}^2}\,
    I_{\zeta\varepsilon}\ecan_{n+1}I_{\eta\varepsilon} =
  \sum_{\varepsilon\in \Edges_{\zeta}} \ \
  \sum_{\buildrel {\scriptstyle\g\in\Omega_{n+1}} \over 
        {\vrule height 7pt width 0pt \g_n=\varepsilon}} 
    {\n{r(\varepsilon)} \over \n{r(\zeta)}}\,
    I_{\zeta\varepsilon}\ecan_{n+1}I_{\eta\varepsilon} \$=
  \sum_{\varepsilon\in \Edges_{\zeta}} \ \
    \n{r(\varepsilon)}\,
    I_{\zeta\varepsilon}\ecan_{n+1}I_{\eta\varepsilon} =
  \sum_{\varepsilon\in \Edges_{\zeta}} \ \
    e^{n+1}_{\zeta\varepsilon,\eta\varepsilon}.
  $$
  Summarizing we have
  $$
  e^n_{\zeta,\eta}=
  \sum_{\varepsilon\in \Edges_{\zeta}} \ \
    e^{n+1}_{\zeta\varepsilon,\eta\varepsilon},
  $$
  from which one easily deduces the proof of the following:

  \state Theorem
  For each $n\in\N$, let $A_n$ be the closed *-sub-algebra of\/ $\CRE$
generated by the set
  $\{e^n_{\g,\d}: \g,\d\in\Omega_n\}$ (see \lcite{\SumMatrix}). Then
$A_n\subseteq A_{n+1}$ and the inclusion of these algebras is
determined by the $n.^{th}$ stage in the Bratteli diagram $\D$ as in
\cite{\Bratteli}.  Therefore the closure of the union of  the
$A_n$ is isomorphic to the AF-algebra associated to $\D$.

In order to prove that $\CRE$ in fact coincides with
$\overline{\bigcup_{n\in\N}A_n}$ it suffices to show that $\CRE$ is
generated by all the 
  $e^n_{\g,\d}$, which we now set out to do.  

\state Proposition The sub-C*-algebra of $\CRE$  generated by the set 
  $$
  \{e^n_{\g,\d}: n\in\N,\ \g,\d\in\Omega_n\}
  $$
  coincides with $\CRE$.

\proof Let $A$ be the sub-C*-algebra of $\CRE$ generated by the set in
the statement.  It is enough to show that $A$ contains $\C$ and all
the $\ecan_n$.  By \lcite{\IandE} we have that every $I_\g\in A$, and
since the set $\{I_\g: n\in\N,\ \g\in\Omega_n\}$ generates $\C$ by the
Stone--Weierstrass Theorem, we have that $\C\subseteq A$.
Given $n\in\N$ notice that 
  $
  1=\sum_{\g\in\Omega_n} I_\g.
  $
  So
  $$
  \ecan_n =
  \Big(\sum_{\g\in\Omega_n} I_\g\Big)\ecan_n\Big(\sum_{\d\in\Omega_n} I_\d\Big) =
  \sum_{\g,\d\in\Omega_n} I_\g\ecan_nI_\d =
  \sum_{\g,\d\in\Omega_n} \n{r(\g)}\inv  e^n_{\g,\d}.
  \proofend
  $$

The proof of \lcite{\Main} therefore follows from our last two
results.

  \section{The groupoid point of view}
  In order to illustrate our description of an AF-algebra as the
C*-algebra of the tail-equivalence relation on its Bratteli diagram we
would now like to show how do our description relates to the
description of AF-algebras as groupoid C*-algebras.

Given a Bratteli diagram $\D$ as above recall that the AF-groupoid
$\G$ associated to $\D$ (see \scite{\Renault}{III.1.1}) coincides with
$\R$, as a set, and the groupoid operations are given as follows:  the
multiplication operation is defined by
  $$
  (\a_1,\a_2) (\b_1,\b_2) = (\a_1,\b_2)
  $$
  whenever $(\a_1,\a_2),(\b_1,\b_2)\in\R$ are such that 
$\a_2 = \b_1$, and the inversion operation is given by
  $$
  (\a,\b)\inv = (\b,\a),
  $$
  for  $(\a,\b)\in\R$.

  The topology of $\G$ is defined as follows: give each $\R_n$ the
product topology (as a subspace of $\Omega\times\Omega$) and say that
a subset $U\subseteq \R$ is open if and only if $U\cap \R_n$ is open
for every $n$.  This defines a topology on $\R$ called the
\stress{inductive limit topology}.

Observe that $\R_n$ is an open subset of $\R_{n+1}$ for every $n$.  In
fact, if $(\a,\b)\in\R_n$ consider the set
  $$
  U =
  \big\{
  (\a',\b') \in \R_{n+1} : \a'_n = \a_n = \b_n=\b'_n
  \big\}.
  $$
  It is easy to see that
  $(\a,\b)\in U\subseteq \R_n$ and since $U$ is clearly open in
$\R_{n+1}$ we conclude that $\R_n$ is indeed open in $\R_{n+1}$.  It
clearly follows that $\R_n$ is open in $\R_m$ for every $m>n$ and
hence we see that
  $\R_n$ is open\fn{By contrast notice that $\R_n$ is not necessarily open in $\R$\/ if the latter is equipped with the product topology.
This is partly the reason why the product topology is not appropriate
and is replaced by the inductive limit topology here.} in $\R$ by
definition of the inductive limit topology.

Observe that the fact that the $\R_n$ are open in $\R$ also ensures
that the inclusion maps
  $$
  \R_n \hookrightarrow \R
  $$
  are homeomorphism onto their images and hence we may view the $\R_n$
(with the product topology) as topological subspaces of $\R$.

Since
  $
  \R_n = \big\{(\a,\b)\in\R_{n+1}: \a_n=\b_n\big\}
  $
  we see that $\R_n$ is closed in $\R_{n+1}$, and hence also in $\R_m$
for every $m>n$.  As before we conclude that $\R_m$ is closed in $\R$.

In particular $\R_0$, the unit space of $\G$, is open in $\R$ and
hence $\G$ is an $r$-discrete groupoid \scite{\Renault}{I.2.6}.

\state Proposition The real valued function $\eb_n$ defined on $\R$ by
  $$
  \eb_n(\a,\b) =
  \left\{\matrix{
    \ds {1 \over \n{s(\a_n)}}, &\hbox{if } (\a,\b)\in\R_n, \cr \cr
    0, &\hbox{otherwise}, }\right.
  $$
  is continuous and of compact support.

\proof Since $\R_n$ is clopen it suffices to prove that $\eb_n$ is
continuous on $\R_n$.  But on $\R_n$ one may describe $\eb_n$ as the
following composition of continuous functions:
  $$
  (\a,\b) \in \R_n \ \longmapsto \a \in \Omega
  \ \longmapsto\ {1 \over \n{s(\a_n)}} \in\Real.
  $$
  Since $\R_n$ is compact it is clear that $\eb_n$ is of compact
support. \proofend 

It is well known that $\G$ is an AF-groupoid which is therefore
\scite{\Renault}{III.1.2} amenable and admits the counting measure as
Haar system.

We will denote by $\CG$ the associated groupoid C*-algebra which is
well known \cite{\Renault} to be isomorphic to the AF-algebra associated to the
Bratteli diagram $\D$.

Recall that the unit space $\R_0$ of $\G$ is a clopen set.  It is also
clear that $\R_0$ is isomorphic to $\Omega$ under the correspondence
  $(\a,\a)\in\R_0\mapsto\a\in\Omega$.
  We will therefore view the algebra $\C$ as a subalgebra of
$\CG$ identifying a function $f\in \C$ with the
continuous function defined on $\G$ by
  $$
  f(\a,\b)=\left\{\matrix{f(\a), & \hbox{ if } \a=\b, \cr\cr
                          0, & \hbox{ if } \a\neq\b.}\right.
  $$

\state Proposition Viewed as elements of $\CG$ the $\eb_n$ above
are projections (self-adjoint idempotent elements) such that
  \izitem
  \zitem $\eb_n\eb_{n+1}=\eb_{n+1}$, and
  \zitem $\eb_n f \eb_n = E_n(f) \eb_n,$
  for every $n$, and every $f\in \C$.

\proof Left to the reader.  \proofend

Since it is clear that $\eb_0$ is
the identity in $\CG$ it follows from
the universal property of $\Toep$ that there exists a unique *-homomorphism 
  $$
  \phi : \Toep \to \CG,
  \eqno{(\seqnumbering)}
  \label\DefinePhi
  $$
  such that $\phi(f)=f$, for all $f\in \C$, and $\phi(\et_n)=\eb_n$,
for all $n$.

For each $n\in\N$, let $\Kb_n$ denote the closure of the image of $\Kt_n$
under $\phi$.  Clearly $\Kb_n$ is the closed linear span of the set
  $$
  \{f\eb_n g: f,g\in\C\}.
  $$
  Given an element $f\eb_n g$ as above notice that, for $(\a,\b)\in\R$,
we have
  $$
  (f\eb_n g)(\a,\b) =
  f(\a)\,\eb_n(\a,\b)\,g(\b),
  $$
  from where we deduce that $f\eb_n g$, if viewed as a continuous
function on $\G$, is supported in $\R_n$.  

Regarding the map
  $
  j: \CG \to C_0(\G)
  $
  described in \scite{\Renault}{II.4.2} we therefore conclude that:

\state \label \SupportInRn Proposition For each $k\in\Kb_n$ one has
that $j(k)$ is supported in $\R_n$.

As a consequence we have:

\state Proposition If $(k_0,\ldots,k_n)$ is a redundancy then
$\phi\(\sum_{i=0}^n k_i\)=0$, where $\phi$ is the map introduced in
\lcite{\DefinePhi}.

\proof Once more referring to the map $j$ of \scite{\Renault}{II.4.2} let 
  $g:= j\(\phi(\sum_{i=0}^n k_i)\)$.  By \lcite{\SupportInRn} we have that $g$
  is supported in $\R_n$.  Since $j$ is one-to-one by
\scite{\Renault}{II.4.2.i} it suffices to show that $g(\a,\b)=0$
whenever $(\a,\b)\in\R_n$.  
Fixing $(\a,\b)\in\R_n$ let $f\in\C$
be such that $f(\b)=1$, while $f(\g)=0$ for all
$\g \in \R_n(\b)\setminus\{\b\}$.

By hypothesis $\(\sum_{i=0}^n k_i\)f\et_n=0$.  Applying  $j\phi$ to
this we conclude that
  $$
  g \star f \star \eb_n=0,
  $$
  where ``$\star$'' refers to the product of
\scite{\Renault}{II.4.2.iii}.  It follows that
  $$
  0 =
  (g \star f \star \eb_n)(\a,\b) =
  \sum_{\g\in\R_n(\b)} g(\a,\g) f(\g) \eb_n(\g,\b) =
  g(\a,\b) \eb_n(\b,\b).
  $$
  Since $\eb_n(\b,\b)\neq0$ we obtain $g(\a,\b)=0$, concluding the
proof.  
  \proofend

\vfill\eject 
It follows that the map $\phi$ of \lcite{\DefinePhi} vanishes on the
redundancy ideal and hence defines a *-homomorphism 
  $$
  \psi : \CRE \to \CG.
  $$

At this point it would not be too difficult to show that $\psi$ is an
isomorphism, e.g.~tracking what happens to the matrix units introduced
above.  After all we now have a very concrete and explicit description
of both $\CRE$ and $\CG$.  We leave the details for the interested
reader.

  \endgroup

  \bye